\documentclass{article}
\usepackage{amsmath,amssymb,amsthm}

\newtheorem{theorem}{Theorem}[section]
\newtheorem{definition}[theorem]{Definition}
\newtheorem{proposition}[theorem]{Proposition}
\newtheorem{corollary}[theorem]{Corollary}

\newcommand{\Nat}{{\mathbb N}}
\newcommand{\Real}{{\mathbb R}}
\newcommand{\Q}{{\mathbb Q}}

\begin{document}
\title{Proof mining in $\Real$-trees and hyperbolic spaces}

\author{Lauren\c tiu Leu\c stean\\[0.2cm] 
\footnotesize Department of Mathematics, DarmstadtUniversity of Technology,\\
\footnotesize Schlossgartenstrasse 7, 64289 Darmstadt, Germany\\[0.1cm]
\footnotesize and\\
\footnotesize Institute of Mathematics "Simion Stoilow'' of the Romanian Academy, \\
\footnotesize Calea Grivi\c tei 21, P.O. Box 1-462, Bucharest, Romania\\[0.1cm]
\footnotesize E-mail: leustean@mathematik.tu-darmstadt.de}
\date{}
\maketitle

\begin{abstract}
This paper is part of the general project of proof mining, developed by Kohlenbach. By "proof mining" we mean the logical analysis of mathematical proofs with the aim of extracting new numerically relevant information hidden in the proofs.

We present logical metatheorems for classes of spaces from functional analysis and hyperbolic geometry, like Gromov hyperbolic spaces, $\Real$-trees and uniformly convex hyperbolic spaces. Our theorems are adaptations to these structures of previous metatheorems of Gerhardy and Kohlenbach, and they guarantee a-priori, under very general logical conditions, the existence of uniform bounds.

We give also an application in nonlinear functional analysis, more specifically in metric fixed-point theory. Thus, we show that the uniform bound on the rate of asymptotic regularity for the Krasnoselski-Mann iterations of nonexpansive mappings in uniformly convex hyperbolic spaces obtained in a previous paper is  an instance of one of our metatheorems.
\end{abstract}
\begin{tabular}{ll}
\footnotesize\noindent {\it Keywords}: &\footnotesize Proof mining, hyeprbolic spaces, $\Real$-trees, asymptotic regularity,\\
&\footnotesize nonexpansive functions
\\
\noindent {\it MSC:\ } & \footnotesize 03F10, 47H09, 20F65, 05C05  
\end{tabular}

\section{Introduction}

This paper is part of the general project of proof mining, developed by Kohlenbach (see \cite{Kohlenbach-book} for details). 

In \cite{Kohlenbach-05-1}, Kohlenbach proved general logical metatheorems which guarantee a-priori, under very general logical conditions, the extractability  of uniform bounds from large classes of proofs in functional analysis, and moreover they provide algorithms for actually extracting effective uniform bounds and transforming the original proof into one for the stronger uniformity results. These metatheorems treat classes of spaces such as metric, hyperbolic spaces in the sense of Reich/Kirk/Kohlenbach (called $W$-hyperbolic spaces in this paper), CAT(0), (uniformly convex) normed, and inner product spaces. They assume the global boundedness of the underlying metric space.  

These metatheorems were vastly generalized in \cite{Gerhardy+Kohlenbach-06} by replacing the assumption of the whole space being bounded with very limited local boundedness assumptions. The new metatheorems guarantee bounds which are uniform for all parameters satisfying these weak local boundedness conditions.   
The proofs are based on a combination between G\" odel's functional interpretation and $a$-majorization, which is a version of majorizability  parametrized by a point $a$ of the space $X$ in question.

In this paper, we present new metatheorems for important structures from hyperbolic geometry or geometric group theory: Gromov hyperbolic spaces, $\Real$-trees and uniformly convex $W$-hyperbolic spaces. These new metatheorems are obtained as adaptations of the existing metatheorems for metric and $W$-hyperbolic spaces, based on the following facts, already noticed in  \cite{Gerhardy+Kohlenbach-06}:
\begin{enumerate}
\item the language may be extended by $a$-majorizable constants. In this case, the extracted bounds then additionally depend on $a$-majorants for the new constants,
\item the theory may be extended by purely universal axioms using new majorizable constants if the types of the quantifiers are appropriate. Then the conclusion holds in all metric spaces $(X,d)$ resp. $W$-hyperbolic spaces $(X,d,W)$ satisfying these axioms (under a suitable interpretation of the new constants if any).
\end{enumerate}

In the last section of the paper we present an application in metric fixed-point theory. In \cite{LL-JMAA-06}, the author obtained a uniform bound on the rate of asymptotic regularity for the Krasnoselski-Mann iterations of nonexpansive mappings in uniformly convex $W$-hyperbolic spaces. These results extend to this more general setting a quantitative version of a strengthening of Groetsch's theorem obtained by Kohlenbach \cite{Kohlenbach-03}. We explain in Section \ref{application}  that the extractability of this uniform bound is an instance of our methateorem for the theory of uniformly convex $W$-hyperbolic spaces. For CAT(0)-spaces, the rate of asymptotic regularity turns out to be quadratic, since  CAT(0)-spaces have a ``nice'' modulus of uniform convexity.

\section{Preliminaries}

In this section we give an informal presentation of the general metatheorem proved by Gerhardy and Kohlenbach \cite{Gerhardy+Kohlenbach-06} for metric spaces and $W$-hyperbolic spaces. We assume familiarity with \cite{Kohlenbach-05-1,Gerhardy+Kohlenbach-06}. The formal system $\mathcal{A}^\omega$ for (weakly extensional) classical analysis is defined in \cite{Kohlenbach-05-1}, and we refer the reader to this paper for all the undefined  notions related to this formal system, as well for the representation of rational and real numbers in $\mathcal{A}^\omega$.

The theory $\mathcal{A}^{\omega}[X,d]_{-b}$ for abstract metric spaces is defined in \cite{Gerhardy+Kohlenbach-06} by extending $\mathcal{A}^{\omega}$ to the set ${\bf T}^{X}$ of all finite types over the ground types $0$ and $X$, adding constants $0_{X}$ of type $X$, and  $d_X$ of type $X\rightarrow X\rightarrow 1$ together with  axioms which make $d_{X}$ a pseudo-metric. Equality $=_{X}$ between objects of type $X$ is defined as:  $x=_{X} y\,:= \, d_{X}(x,y)=_{\Real} 0_{\Real}$.

We present in the sequel the setting of hyperbolic spaces as introduced by Kohlenbach \cite{Kohlenbach-05-1}; see \cite{Kohlenbach-05-1,Kohlenbach-book} for detailed discussion of this and related notions. In order to distinguish them from the usual notion of hyperbolic space from hyperbolic geometry and from Gromov hyperbolic spaces, we shall call them $W$-hyperbolic spaces. 

A {\em  W-hyperbolic space}
is a triple $(X,d,W)$ where
$(X,d)$ is a metric space and $W:X\times X\times [0,1]\to X$ is such that for all $x,y,z,w\in X, \lambda_1,\lambda_2\in [0,1],$
\[\begin{array}{ll}
(W1)\quad & d(z,W(x,y,\lambda))\le (1-\lambda)d(z,x)+\lambda d(z,y),\\
(W2)\quad & d(W(x,y,\lambda_1),W(x,y,\lambda_2))=|\lambda_1-\lambda_2|\cdot d(x,y),\\
(W3)\quad & W(x,y,\lambda)=W(y,x,1-\lambda),\\
(W4)\quad & d(W(x,z,\lambda),W(y,w,\lambda)) \le (1-\lambda)d(x,y)+\lambda d(z,w).
\end{array}\]
If $x,y\in X$ and $\lambda\in[0,1]$, then we use the notation $(1-\lambda)x\oplus \lambda y$ for $W(x,y,\lambda)$. We shall denote by $[x,y]$ the set $\{(1-\lambda)x\oplus \lambda y:\lambda\in[0,1]\}$. A nonempty subset $C\subseteq X$ is {\em convex} if $[x,y]\subseteq C$ for all $x,y\in C$. 

The theory $\mathcal{A}^\omega[X,d,W]_{-b}$ results from $\mathcal{A}^\omega[X,d]_{-b}$ by adding  a new constant $W_X$ of type $X\rightarrow X\rightarrow 1 \rightarrow X$ together with the appropriate axioms.

If $X$ is a nonempty set, the full-theoretic type structure $S^{\omega,X}:=\left<S_\rho\right>_{\rho\in {\bf T}^X}$ over $\Nat$ and $X$ is defined by $S_0:=\Nat, \quad S_X:=X, S_{\rho\rightarrow \tau}:=S_\tau^{S_\rho},$
where $S_\tau^{S_\rho}$ is the set of all set-theoretic functions $S_\rho\rightarrow S_\tau$.
\begin{definition}\label{model-metric+W-hyp}
We say that a sentence of ${\cal L}(\mathcal{A}^\omega[X,d]_{-b})$ holds in a nonempty metric space $(X,d)$ if it holds in the models of $\mathcal{A}^\omega[X,d]_{-b}$ obtained as follows: by letting the variables range over the appropriate universe of the full-theoretic type structure $S^{\omega,X}$ with the set $X$ as the universe for the base type $X$; $0_X$ is interpreted by an arbitrary element of $X$; $d_X$ is interpreted as $d_X(x,y):=(d(x,y))_\circ$.

The notion that a sentence of ${\cal L}(\mathcal{A}^\omega[X,d,W]_{-b})$ holds in a nonempty $W$-hyperbolic space $(X,d,W)$ is obtained from the previous one by interpreting  $W_X(x,y,\lambda^1)$ as $W(x,y,r_{\tilde{\lambda}})$, where $r_{\tilde{\lambda}}\in[0,1]$ is the unique real number represented by $\tilde{\lambda}$.
\end{definition}
In the above definition, $\lambda\mapsto \tilde{\lambda}$ is used \cite{Kohlenbach-05-1} to represent the interval $[0,1]$ by number-theoretic functions $\Nat\to\Nat$ and  $(\cdot)_{\circ}$ is a semantic operator, defined also in \cite{Kohlenbach-05-1}, which for any real number $x\in[0,\infty)$ selects out of all the representatives $f\in \Nat^{\Nat}$ of $x$  a unique representative $(x)_{\circ}\in\Nat^{\Nat}$  satisfying some ``nice'' properties: 
\[(x)_{0}(n):=j(2k_{0}, 2^{n+1}-1), \]
where $\displaystyle k_0=\max k\left[\frac k {2^{n+1}}\leq x\right]$,  and $j$ is the Cantor pairing function.

\begin{definition}\cite{Kohlenbach-CIE}
A type $\rho$ is called {\em small} if it is of degree $\leq 1$ (i.e. $0\rightarrow \ldots \rightarrow 0$) or of the form $\rho_1\rightarrow\ldots\rightarrow\rho_k\rightarrow X$ with the $\rho_i$ being of type $0$ or $X$.
\end{definition}

\begin{definition}\cite{Kohlenbach-CIE}
A formula $A$ is called a $\forall$-formula (resp. $\exists$-formula) if it has the form $A\equiv\forall\underline{x}^{\underline{\sigma}}A_0(\underline{x})$ (resp. $A\equiv \exists \underline{x}^{\underline{\sigma}}A_0(\underline{x})$) where $A_0$  is a quantifier free formula and the types in $\underline{\sigma}$ are small.
\end{definition}

For any type $\rho\in {\bf T}^X$, we define the type $\widehat{\rho}\in {\bf T}$, which is the result of replacing all occurrences of the type $X$ in $\rho$ by the type $0$. Based on Bezem's notion of strong majorizability {\em s-maj} \cite{Bezem}, Gerhardy and Kohlenbach \cite{Gerhardy+Kohlenbach-06} defined a parametrized $a$-majorization relation $\gtrsim^a$ between objects of type $\rho\in {\bf T}^X$ and their majorants of type $\widehat{\rho}\in {\bf T}$, where the parameter $a$ of type $X$ serves as a reference point for comparing and majorizing elements of $X$:
\[ {x^*}^0 \gtrsim^a_0  x^0:\equiv  x^* \ge_0 x, \quad\quad {x^*}^0\gtrsim^a_X x^X  :\equiv (x^*)_\Real\ge_\Real d_X(x,a),\]
\[x^* \gtrsim^a_{\rho\rightarrow\tau} x :\equiv  \forall y^*,y(y^*\gtrsim^a_{\rho} y \rightarrow x^* y^* \gtrsim^a_{\tau} xy)\wedge \forall z^*,z(z^*\gtrsim^a_{\hat{\rho}} z \rightarrow x^* z^* \gtrsim^a_{\hat{\tau}} x^*z).
\]
Restricted to the types $\bf T$ the relation $\gtrsim^a$ is identical with the strong majorizability {\em s-maj} and, hence, for $\rho\in \bf T$ we shall write {\em s-maj}$_\rho$ instead of $\gtrsim^a_\rho$, since in this case the parameter $a$ is irrelevant.

The following theorem is a simplified version of the very general metatheorem proved first by Kohlenbach \cite{Kohlenbach-05-1} for bounded metric ($W$-hyperbolic) spaces, and then generalized to the unbounded case  by Gerhardy and Kohlenbach: 

\begin{theorem}\label{main-thm}\cite{Gerhardy+Kohlenbach-06}
Let $\rho$ be a small type and  $B_\forall(x^\rho,n^0)$ (resp. $C_\exists(x^\rho,m^0)$) be a $\forall$-formula containing  only $x,n$ free (resp. a $\exists$-formula containing only $x,m$ free). Assume that the constant $0_X$ does not occur in $B_\forall, C_\exists$ and that
\begin{equation}
\mathcal{A}^\omega[X,d]_{-b}\vdash \forall x^\rho (\forall n B_\forall(x,n)\rightarrow \exists m C_\exists(x,m)).
\end{equation}
Then there exists a computable functional $\Phi:S_{\widehat{\rho}}\rightarrow \Nat$ such that the following holds in all nonempty metric spaces $(X,d)$:

for all $x\in S_\rho, x^*\in S_{\widehat{\rho}}$, if there exists an $a\in X$ such that $\displaystyle x^{*}\gtrsim^a x$, then
\begin{equation}
\forall n\le \Phi(x^*) B_\forall(x,n)\rightarrow \exists m\le \Phi(x^*)C_\exists(x,m).
\end{equation}
The theorem also holds for $\mathcal{A}^\omega[X,d,W]_{-b}$ and nonempty $W$-hyperbolic spaces $(X,d,W)$.
\end{theorem}
Instead of single variables $x,n,m$ and single premises $\forall n B_\forall(x,n)$ we may have tuples of variables and finite conjunctions of premises. In the case of a tuple $\underline{x}$, we have to require that we have a tuple $\underline{x}^*$ of $a$-majorants for a common $a\in X$ for all the components of the tuple $\underline{x}$.

\section{New metatheorems}

In this section, we extend Theorem \ref{main-thm} to important structures from hyperbolic geometry or geometric group theory: Gromov hyperbolic spaces, $\Real$-trees and uniformly convex $W$-hyperbolic spaces. 

As we have already remarked, in order to get metatheorems for these new structures, it is enough to reformulate their theories as extensions of $\mathcal{A}^\omega[X,d]_{-b}$ or $\mathcal{A}^\omega[X,d,W]_{-b}$ by purely universal sentences having appropriate types for quantifiers and to verify that the new constants (if any) are $a$-majorizable. We shall see in the sequel that this is possible.

\subsection{Gromov hyperbolic spaces}

Gromov's theory of hyperbolic spaces is set out in \cite{Gromov-87}.  The study of Gromov hyperbolic spaces has been largely motivated and dominated by questions about (Gromov) hyperbolic groups, one of the main object of study in geometric group theory. In the sequel, we review some definitions and elementary facts concerning Gromov hyperbolic spaces. For a more detailed account of this material, the reader is referred to \cite{Gromov-87,Ghys+deLaHarpe-eds,Bridson+Haefliger}. 

Let $(X,d)$ be a metric space. Given three points $x,y,w$, the {\em Gromov product} of $x$ and $y$  with respect to the {\em base point} $w$ is defined to be:

$\quad\quad(x\cdot y)_w=\frac 12(d(x,w)+d(y,w)-d(x,y)).$\\
It measures the failure of the triangle inequality to be an equality and it is always nonnegative. 
\begin{definition}\label{delta-hyp}
Let $\delta\geq 0$. $X$ is called $\delta-hyperbolic$ if for all $x,y,z,w\in X$,
\begin{equation}
(x\cdot y)_w\geq \min\{(x\cdot z)_w, (y\cdot z)_w\}-\delta.\label{delta-hyp-ineq}
\end{equation}
We say that $X$ is {\em Gromov hyperbolic} if  it is $\delta$-hyperbolic
for some $\delta\geq 0$.
\end{definition}
It turns out that the definition is independent of the choice of the base point $w$ in the sense that if the Gromov product is $\delta$-hyperbolic with respect to one base point then it is $2\delta$-hyperbolic with respect to any base point.

By unraveling the definition of Gromov product, (\ref{delta-hyp-ineq}) can be rewritten as a 4-point condition: for all $x,y,z,w\in X$,
\begin{equation}
d(x,y)+d(z,w)\leq \max\{d(x,z)+d(y,w), d(x,w)+d(y,z)\}+2\delta. \label{Q(delta)-ineq}
\end{equation}

The theory of Gromov hyperbolic spaces, $\mathcal{A}^\omega[X,d,\delta\text{-hyperbolic}]_{-b}$ is defined by extending  $\mathcal{A}^\omega[X,d]_{-b}$ as follows: 
\begin{enumerate}
\item add a constant $\delta^1$ of type $1$,
\item add the axioms 

$\begin{array}{l}
\delta \ge_\Real 0_\Real,\\
\forall x^X, y^X, z^X, w^X\,\big(d_X(x,y)+_\Real d_X(z,w)\leq_{\Real} \\
 \quad  \quad \leq_\Real{\rm max}_\Real\{d_X(x,z)+_\Real d_X(y,w),d_X(x,w)+_\Real d_X(y,z)\}+_\Real 2\cdot_\Real\delta\big). 
\end{array}$
\end{enumerate}
The notion that a sentence of ${\cal L}(\mathcal{A}^\omega[X,d,\delta\text{-hyperbolic}]_{-b})$ holds in a nonempty  Gromov hyperbolic space $(X,d,\delta)$ is defined as in Definition \ref{model-metric+W-hyp}, by interpreting the new constant $\delta^1$ as $\delta^1:=(\delta)_0$.

Since $\leq_\Real$ is $\Pi_0^1$, the two axioms are universal. Thus, in order to adapt Theorem \ref{main-thm} to the theory of Gromov hyperbolic spaces, we need to show that the new constant $\delta^1$ is strongly majorizable. It is easy to see that if $(X,d,\delta)$ is a $\delta$-hyperbolic space, and $k\in\Nat$ is such that $k\ge \delta$, then 
\[\delta^*:=\lambda n.j(k\cdot 2^{n+2}, 2^{n+1}-1) \,\, \text{s-maj}_1 (\delta)_\circ.\] 

\begin{theorem}\label{main-thm-delta-hyp}
Theorem \ref{main-thm} holds also for $\mathcal{A}^\omega[X,d, \delta\text{-hyperbolic}]_{-b}$ and non-empty  Gromov hyperbolic spaces $(X,d,\delta)$, with the bound $\Phi$ depending additionally on  $k\in\Nat$ such that $k\geq\delta$.
\end{theorem}

\subsection{$\Real$-trees}

The notion of $\Real$-tree was introduced by Tits \cite{Tits}, as a generalization of the notion of local Bruhat-Tits building for rank-one groups, which itself generalizes the notion of simplicial tree. A more general concept, that of a $\Lambda$-tree, where $\Lambda$ is a totally ordered abelian group,  made its appearance as an essential tool in the study of groups acting on hyperbolic manifolds in the work of Morgan and Shalen \cite{Morgan+Shalen-84}. For detailed informations about $\Real(\Lambda)$-trees, we refer to \cite{Chiswell-01}.

In the sequel we recall some basic definitions. Let $(X,d)$ be a metric space. A {\em geodesic} in $X$ is a map $\gamma:[a,b]\to\Real$ which is distance-preserving, that is $ d(\gamma(s),\gamma(t))=|s-t| \text{~~for all~~} s,t\in [a,b].$ A {\em geodesic segment} in $X$ is the image of a geodesic in $X$. 
If  $\gamma :[a,b]\to\Real$ is a geodesic, and $x,y\in X$ are such that $\gamma(a)=x$ and $\gamma(b)=y$, we say that $\gamma$ is a {\em geodesic from  x to y} or that the geodesic segment $\gamma([a,b])$ {\em joins x and y}.  $X$ is said to be a  {\em (uniquely) geodesic space} if every two points are joined by a (unique) geodesic. If $X$ is a uniquely geodesic space, then we denote by $[x,y]$ the unique geodesic segment that joins $x$ and $y$.

\begin{definition}\cite{Tits}
$X$ is an $\Real$-tree iff  $X$ is a  geodesic space containing no homeomorphic image of a circle.
\end{definition}

We remark that in the initial definition, Tits only considered $\Real$-trees which are complete as metric spaces, but the assumption of completeness is usually irrelevant. The following proposition gives some equivalent characterizations of $\Real$-trees, which can be found in the literature.

\begin{proposition}\label{R-trees-equiv-char}
Let $(X,d)$ be a metric space. The following are equivalent:
\begin{enumerate}
\item $X$ is an $\Real$-tree,
\item $X$ is uniquely geodesic and for all $x,y,z\in X$, 
\begin{equation} 
[y,x]\cap[x,z]=\{x\}\Rightarrow [y,x]\cup[x,z]=[y,z].
\end{equation}
(i.e., if two geodesic segments intersect in a single point, then their union is a geodesic segment.)
\item $X$ is a geodesic space which is $0$-hyperbolic, i.e. satisfies the inequality (\ref{Q(delta)-ineq}) with $\delta=0$.
\end{enumerate}
\end{proposition}

The fact that $\Real$-trees are exactly the geodesic $0$-hyperbolic spaces follows from a very important result of Alperin and Bass \cite[Theorem 3.17]{Alperin+Bass-87} and is the basic ingredient for proving the following characterization of $\Real$-trees using our notion of $W$-hyperbolic space.
 
\begin{proposition}
Let $(X,d)$ be a metric space. The following are equivalent:
\begin{enumerate}
\item $X$ is an $\Real$-tree;
\item $X$ is a $W$-hyperbolic space which satisfies for all $x,y,z,w\in X$,
\begin{eqnarray*}
d(x,y)+d(z,w)\leq \max\{d(x,z)+d(y,w),d(x,w)+d(y,z)\}.
\end {eqnarray*}
\end{enumerate}
\end{proposition}

Now, we are ready to define the formal theory $\mathcal{A}^\omega[X,d,W,\Real\text{-tree}]_{-b}$ of $\Real$-trees. This results from the theory $\mathcal{A}^\omega[X,d,W]_{-b}$ by adding the axiom:
\[
\,\,\begin{cases}\forall x^X,y^X,z^X,w^X\big(d_X(x,y)+_\Real d_X(z,w)\leq_\Real \\
\quad\quad\quad\quad\quad\quad\leq_\Real {\rm max}_\Real\{d_X(x,z)+_\Real d_X(y,w),d_X(x,w)+_\Real d_X(y,z)\}\big).
\end{cases}\]

Hence, $\mathcal{A}^\omega[X,d,W,R\text{-tree}]_{-b}$  is obtained from $\mathcal{A}^\omega[X,d,W]_{-b}$ only by adding an universal axiom.

\begin{theorem}\label{main-thm-R-tree}
Theorem \ref{main-thm} holds also for $\mathcal{A}^\omega[X,d,W,R\text{-tree}]_{-b}$ and nonempty $\Real$-trees. 
\end{theorem}

\subsection{Uniformly convex W-hyperbolic spaces}

The notion of uniformly convex $W$-hyperbolic space is defined in \cite{LL-JMAA-06}, following \cite[p.105]{Goebel/Reich(84)}. 

\begin{definition}
A $W$-hyperbolic space $(X,d,W)$ is called {\em uniformly convex} if for 
any $r>0$, and $\varepsilon\in(0,2]$ there exists a $\delta\in(0,1]$ such that 
for all $a,x,y\in X$,
\begin{eqnarray}
\left.\begin{array}{l}
d(x,a)\le r\\
d(y,a)\le r\\
d(x,y)\ge\varepsilon r
\end{array}
\right\}
& \quad \Rightarrow & \quad d\left(\frac12x\oplus\frac12y,a\right)\le (1-\delta)r. \label{uc-def}
\end{eqnarray}
A mapping $\eta:(0,\infty)\times(0,2]\rightarrow (0,1]$ providing such a
$\delta:=\eta(r,\varepsilon)$ for given $r>0$ and $\varepsilon\in(0,2]$ is called a {\em modulus of uniform
convexity}.
\end{definition}

Using standard continuity arguments, we can prove the following equivalent characterization.

\begin{proposition}
Let $(X,d,W)$ be a $W$-hyperbolic space. The following are equivalent:
\begin{enumerate}
\item there exists $\eta:(0,\infty)\times(0,2]\rightarrow (0,1]$ such that (\ref{uc-def}) holds,
\item there exists $\eta:\Q_*^+\times \Nat\rightarrow \Nat$ such that 
for any $r\in\Q_*^+,k\in\Nat$, and $a,x,y\in X$ 
\begin{eqnarray}
\left.\begin{array}{l}
d(x,a)< r\\
d(y,a) < r\\
d\left(\frac12x\oplus\frac12y,a\right) > \left(1-2^{-\eta(r,k)}\right)r
\end{array}
\right\}
& \quad \Rightarrow & \quad d(x,y)\le 2^{-k}r. 
\end{eqnarray}
\end{enumerate}
\end{proposition}

The theory $\mathcal{A}^\omega[X,d,W,\eta]_{-b}$ of uniformly convex $W$-hyperbolic spaces extends the theory $\mathcal{A}^\omega[X,d,W]_{-b}$ as follows:
\begin{enumerate}
\item add a new constant $\eta_X$ of type $000$,
\item add the following axioms:\\
$\begin{array}{l}
\forall r^0\forall k^0\forall x^X,y^X,a^X \big(d_X(x,a) <_\Real r \,\wedge\, d_X(y,a)<_\Real r\,\wedge\,\\
\wedge\, \displaystyle d_X(W_X(x,y,1/2),a)>_\Real \displaystyle \left(1-2^{-\eta_X(r,k)}\right)\cdot_\Real r\rightarrow d_X(x,y)\leq_\Real 2^{-k}\cdot_\Real r\big), \\
\forall r^0,k^0(\eta_X(r,k)=_0\eta_X(q(r),k)). 
\end{array}$
\end{enumerate}
In the second axiom, we express the fact that $\eta_X$ is a function having the first argument a rational number on the level of codes. The function $q$ is defined by: $q(n):=\min k\leq_0 n[k=_\Q n]$ (see \cite{Kohlenbach-entcs-98} for details).
Since $<_\Real\in \Sigma_0^1$ and $\leq_\Real\in \Pi_0^1$, it follows that our new theory is obtained from $\mathcal{A}^\omega[X,d,W]_{-b}$ by adding two universal axioms. It is easy to see also that  the constant $\eta_X^{000}$ is majorizable. 

The notion that a sentence of ${\cal L}(\mathcal{A}^\omega[X,d,W,\eta]_{-b})$ holds in a nonempty uniformly convex $W$-hyperbolic space $(X,d,W,\eta)$ is defined as in Definition \ref{model-metric+W-hyp}, by interpreting the new constant $\eta_X$ as $\eta_X(r,k):=\eta(q(r),k)$. 

\begin{theorem}\label{main-thm-uc-hyp}
Theorem \ref{main-thm} holds also for $\mathcal{A}^\omega[X,d,W,\eta]_{-b}$ and nonempty uniformly convex $W$-hyperbolic spaces $(X,d,W,\eta)$, with the bound $\Phi$ depending additionally on the modulus of uniform convexity $\eta$.
\end{theorem}

\section{An application to metric fixed point theory}\label{application}

In this last section we present an application of the metatheorem for uniformly convex $W$-hyperbolic spaces to metric fixed point theory, more specifically to fixed point theory of nonexpansive functions.

Let $(X,d,W)$ be a hyperbolic space and $C\subseteq X$ a nonempty convex subset of $X$. 
A mapping $T:C\to C$ is called {\em nonexpansive} ({\em n.e.} for short) if $d(Tx,Ty)\le d(x,y)$ for all $x,y\in C$. As in the case of normed spaces \cite{Mann(53),Krasnoselski(55)},  we can define the {\em Krasnoselski-Mann iteration} starting from $x\in C$ by: 
\begin{equation}
x_0:=x, \quad x_{n+1}:=(1-\lambda_n)x_n \oplus\lambda_n Tx_n, \label{KM-lambda-n-def-hyp}\end{equation}
where $(\lambda_n)$ is a sequence in $[0,1]$.

Asymptotic regularity  was defined by Browder and Petryshyn \cite{Browder/Petryshyn(66)}:  a mapping $T:C\to C$ is called {\em asymptotically regular}  if  
$\displaystyle\lim_{n\to\infty}d(T^n(x),T^{n+1}(x))=0$ for all $x\in C$. Following \cite{Borwein/Reich/Shafrir(92)}, we say that a nonexpansive mapping $T:C\to C$ is {\em $\lambda_n$-asymptotically regular} if  $\displaystyle\lim_{n\to\infty}d(x_n,Tx_n)=0$ for all $x\in C$. The following theorem was proved by the author in \cite{LL-JMAA-06}:

\begin{theorem}\label{Groetsch-hyp}
Let $(X,d,W)$ be a uniformly convex hyperbolic space  with modulus of uniform convexity $\eta$ such that $\eta$ decreases with  $r$ (for a fixed $\varepsilon$), $C\subseteq X$ be a nonempty convex subset and $T:C\rightarrow C$ nonexpansive such that $T$ has at least one fixed point. 
Assume moreover that $(\lambda_n)$ is a sequence in $[0,1]$ such that $\displaystyle\sum_{n=0}^\infty\lambda_n(1-\lambda_n)=\infty. $

Then $T$ is $\lambda_n$-asymptotically regular.
\end{theorem}

This theorem is the version for uniformly convex $W$-hyperbolic spaces of a theorem for normed spaces proved by Groetsch \cite{Groetsch(72)}. In Groetsch's theorem the extra-hypothesis on $\eta$ is not needed, due to nice scaling properties of normed spaces.  Still, we must emphasize that this extra-hypothesis is satisfied by important classes of  uniformly convex $W$-hyperbolic spaces like the Hilbert ball, CAT(0)-spaces or $\Real$-trees.

In the sequel, we show that Theorem \ref{main-thm-uc-hyp}  guarantees uniform effective bounds for a strengthening of the above theorem, which only assumes the existence of approximate fixed points in some neighborhood of the starting point $x$.

Since any convex subset of a $W$-hyperbolic space is also a $W$-hyperbolic space, it suffices to consider only nonexpansive functions $T:X\to X$.

We use the following notations:
\[\begin{array}{lll}
Mon(\eta,r)&:= &\forall r_1^0,r_2^0,k^0(r_1\leq_\Q r_2\rightarrow \eta(r_1,k)\geq_0 \eta(r_2,k)),\\
Fix(T)&:= &\{p^X \mid T(p)=_X p\},\\
Fix_\delta(T,x,b)&:= &\{y^X\mid d_X(y,T(y))\le_\Real \delta \wedge d_X(x,y)\le_\Real b\}.
\end{array}\]

The following more concrete consequence of Theorem \ref{main-thm-uc-hyp} suffices for our application. Its proof is similar with the one of \cite[Corollary 4.22]{Gerhardy+Kohlenbach-06}.
 
\begin{corollary}\label{cor-application}
Let $P$ be a $\mathcal{A}^\omega$-definable Polish space and $K$ be a $\mathcal{A}^\omega$-definable compact Polish space. Let $B_\forall$ and $C_\exists$ be as before. If
$\mathcal{A}^\omega[X,d,W,\eta]_{-b}$ proves that
\begin{eqnarray*}
\forall z\in P\forall y\in K\forall x^X, T^{X\to X}\big(T \,n.e.  \wedge Fix(T)\ne \emptyset \wedge \forall n^0B_\forall \rightarrow \exists m^0C_\exists),
\end {eqnarray*}
then there exists a computable functional $\Phi:\Nat^\Nat\times\Nat\times \Nat^{\Nat\times \Nat}\to\Nat $ (on representatives $r_z:\Nat\to\Nat$ of elements $z\in P$) such that for all $r_z\in\Nat^\Nat, b\in\Nat$
\begin{eqnarray*}
\forall y\in K\forall x^X, T^{X\to X}\big(T\, n.e.\, \wedge\, \forall \delta>0\left(Fix_\delta(T,x,b)\ne\emptyset\right)\,\wedge\\
\wedge\,\forall n\leq_0 \Phi(r_z,b,\eta)B_\forall \rightarrow \exists m\leq_0 \Phi(r_z,b,\eta)C_\exists)
\end {eqnarray*}
holds in any nonempty uniformly convex $W$-hyperbolic space $(X,d,W,\eta)$.

As before, instead of single variables $y,z$ and single premises $B_\forall$, we may have tuples of variables and a finite conjunction of premises. 
\end{corollary}

Using the fact that the sequence $(d(x_n,T(x_n))$ is nonincreasing, we get that
T is asymptotically regular is equivalent with
\begin{eqnarray*}
\forall x\in X\forall k\in\Nat \exists N\in \Nat(d(x_N, T(x_N))< 2^{-k}).
\end {eqnarray*}
The assumption  on $(\lambda_n)$ in Theorem \ref{Groetsch-hyp} is equivalent with the existence of a witness $\theta:\Nat\to\Nat$ such that for all $n\in \Nat$,
$\displaystyle\sum_{i=0}^{\theta(n)} \lambda_i(1-\lambda_i)\geq n.$

It follows that $\mathcal{A}^\omega[X,d,W,\eta]_{-b}$ proves the following formalized version of Theorem \ref{Groetsch-hyp}:
\begin{eqnarray*}
\forall k^0  \,\forall \theta^1 \,\forall\lambda^{0\to 1}_{(\cdot )}\,\forall x^X, T^{X\to X}\,\big(Mon(\eta,r) \,\wedge \, T\, n.e. \,\wedge \, Fix(T)\ne\emptyset\, \wedge\quad\quad\\
\wedge \forall n^0(n\leq_\Real\sum_{i=0}^{\theta(n)}\lambda_i(1-\lambda_i))\rightarrow \exists N^0(d_X(x_N, T(x_N))<_\Real 2^{-k})\big),
\end{eqnarray*}
where $\lambda^{0\to 1}_{(\cdot )}$ represents an element of the compact Polish space $[0,1]^\infty$ with the product metric.  Corollary \ref{cor-application} yields the existence of a computable functional $\Phi(k,\theta, b, \eta)$ such that for all $(\lambda_n)\in[0,1]^\infty, x\in X, T:X\to X$,
\begin{eqnarray*}
\!\!\!\!\!\!\!\!\!\!\!\!Mon(\eta,r)\wedge T\, n.e.\wedge\forall \delta>0(Fix_\delta(T,x,b)\ne\emptyset)\wedge\forall n(n\leq \sum_{k=0}^{\theta(n)} \lambda_k(1-\lambda_k))\rightarrow\\
\rightarrow \exists N\leq \Phi(k,\theta, b, \eta)(d(x_N, T(x_N))\leq  2^{-k})\big)
\end{eqnarray*}
holds  in any nonempty uniformly convex $W$-hyperbolic space $(X,d,W, \eta)$.
Using again that $(d(x_n,T(x_n)))$ is nonincreasing, it follows that $\Phi(k,\theta, b, \eta)$ is a bound on the rate of convergence of $(d(x_n,T(x_n)))$ towards $0$.

Hence, as an application of Corollary \ref{cor-application}, we immediately obtain the following uniform version of a strengthening of Theorem \ref{Groetsch-hyp}.

\begin{theorem}
Let $(X,d,W)$ be a uniformly convex hyperbolic space  with modulus of uniform convexity $\eta$ such that $\eta$ decreases with  $r$ (for a fixed $\varepsilon$), $C\subseteq X$ be a nonempty convex subset and $T:C\rightarrow C$ nonexpansive.  
Assume that $(\lambda_n)$ is a sequence in $[0,1]$ and $\theta :\Nat\to\Nat$ is such that for all $n\in\Nat$, 
\begin{equation}
\sum\limits_{i=0}^{\theta (n)} \lambda_i(1-\lambda_i) \ge n. 
\label{main-hyp-lambda-gamma}
\end{equation} 
Let $x\in C,b>0$ be such that for any $\delta >0$ there is $y\in C$ with 
\begin{equation}
\rho(x,y)\le b \quad \mbox{and~~}  \rho(y,Ty) \leq\delta. \label{main-hyp-x-y}
\end{equation}
Then $\displaystyle\lim_{n\to\infty} \rho(x_n,Tx_n)=0,$ and moreover
\[
\forall \varepsilon >0\, \forall n\ge \Phi(\varepsilon,\theta,b,\eta)\, \big( \rho(x_n,Tx_n)\leq \varepsilon\big). \label{main-thm-conclusion}
\]
\end{theorem}
The extraction of $\Phi(\varepsilon,\theta, b, \eta)$ was carried out in \cite[Theorem 14]{LL-JMAA-06},  :
\[\Phi(\varepsilon,\theta,b,\eta):=\left\{ \begin{array}{ll}
        \displaystyle \theta\left(\left\lceil\frac{b+1}{\varepsilon\cdot\eta\left(b+1,\displaystyle\frac{\varepsilon}{b+1}\right)}
\right\rceil\right) & \text{for~ } \varepsilon <2b\\
        0 & \text{otherwise. }
        \end{array}\right.
\]

Moreover, for bounded $C$, the condition (\ref{main-hyp-x-y}) holds for all $x\in C$ with $d_C$ instead of $b$, so we get asymptotic regularity for general $(\lambda_n)$ satisfying (\ref{main-hyp-lambda-gamma}) and an explicit bound $\Phi(\varepsilon,\theta,d_C,\eta)$ on the rate of asymptotic regularity, which depends only on the error $\varepsilon$, on the modulus of uniform convexity $\eta$,  on the diameter $d_C$ of $C$ and on $(\lambda_n)$ only via $\theta$, but not  on the nonexpansive mapping $T$, the starting point $x\in C$ of the iteration or other data related with $C$ and $X$. Furthermore, for $CAT(0)$-spaces (and subsequently for $\Real$-trees), which have a very "nice'' modulus of uniform convexity, we have got a quadratic  rate of asymptotic regularity. I refer the reader to \cite{LL-JMAA-06} for a detailed presentation of all these facts.\\

\noindent {\bf Acknowledgement}\\

\noindent I express my gratitude to Ulrich Kohlenbach for the numerous discussions we had on the subject. His suggestions led to a substantially improved presentation of the results.

\end{document}